\documentclass[12pt]{article}
\usepackage{amsfonts}
\usepackage{mathrsfs}
\usepackage{amsmath,amssymb}
\numberwithin{equation}{section}

\newtheorem{theo}{Theorem}[section]
\newtheorem{conj}[theo]{Conjecture}
\newtheorem{prob}[theo]{Problem}

\newtheorem{lemm}[theo]{Lemma}

\begin{document}

\centerline{\bf Irreducible weight modules with a finite-dimensional weight space}
\centerline{\bf{over the twisted $N\!=\!1$ Schr\"{o}dinger-Neveu-Schwarz algebra}}

\vspace*{6pt}

\centerline{Huanxia Fa$^{1)}$, Jianzhi Han$^{2)}$, Junbo Li$^{1)}$}

\centerline{\small $^{1)}$School of Mathematics and Statistics, Changshu Institute of Technology, Changshu 215500, China}

\centerline{\small $^{2)}$Department of Mathematics, Tongji University, Shanghai, 200092, China}

\centerline{\small E-mail: sd\_huanxia@163.com, jzhan@tongji.edu.cn, sd\_junbo@163.com}

\vspace*{6pt}

\noindent{\small{\bf Abstract.}
It is shown that there are no simple mixed modules over the twisted $N\!=\!1$ Schr\"{o}dinger-Neveu-Schwarz algebra, which implies that every irreducible weight module over it with a nontrivial finite-dimensional weight space, is a Harish-Chandra module.


}

\vspace*{10pt}

\centerline{\bf1\ \ Introduction}
\setcounter{section}{1}\setcounter{theo}{0}

The Schr\"{o}dinger-Virasoro type Lie algebras and their supersymmetric counterparts are closely related to mathematics and physics. The Schr\"{o}dinger-Virasoro algebra was originally introduced by M. Henkel in 1994 during the process of trying to apply the concepts and methods of conformal field theory to models of statistical physics (see \cite{H--JSM1994} for reference). The deformative counterparts were introduced in \cite{RU--AHP2006} and their supersymmetric extensions were investigated in the context of supersymmetric quantum mechanics. Schr\"{o}dinger-Neveu-Schwarz Lie superalgebras $\frak{sns}^{(N)}$ with $N$ supercharges were introduced in \cite{HU--NPB2006} generally, which appear as a semi-direct product of Lie algebras of super-contact vector fields with infinite-dimensional nilpotent Lie superalgebras. The {\it twisted $N\!=\!1$ Schr\"{o}dinger-Neveu-Schwarz algebra} $\widetilde{\frak{tsns}}$ (introduced firstly in \cite{HU--NPB2006}) is an
infinite-dimensional Lie superalgebra over $\mathbb{C}$ with the
basis $\{L_n,G_r,Y_p,M_p,\frak{c}\,|\,n\in \mathbb{Z},\,r\in\frac{1}{2}+\mathbb{Z},\,p\in\frac{1}{2}\mathbb{Z}\}$ and the following non-vanishing super brackets:
\begin{eqnarray*}
&&[L_n,L_{m}]=(m-n)L_{m+n}+\frac{n^3-n}{12}\delta_{m,-n}\frak{c},\\[6pt] &&[L_n,G_r\,]=(r-\frac{n}{2})G_{r+n},\ \ \ \ [G_r,G_s]=2L_{r+s}+\frac{1-4r^2}{12}\delta_{r,-s}\frak{c},\\[6pt]
&&[L_n,Y_p]=\left\{\begin{array}{cc}
(p-\frac{n}{2})Y_{p+n},&{\rm if}\ p\in\mathbb{Z},\\[6pt]
pY_{p+n},&{\rm if}\ p\in\frac{1}{2}+\mathbb{Z},\end{array}\right.
\ \ \,[L_n,M_p]=\left\{\begin{array}{cc}
pM_{p+n},&{\rm if}\ p\in\mathbb{Z},\\[6pt]
(p+\frac{n}{2})M_{p+n},&{\rm if}\ p\in\frac{1}{2}+\mathbb{Z},\end{array}\right.\\
&&[G_r,Y_p]=
\left\{\begin{array}{cc}
\frac{1}{2}(p-r)Y_{p+r},&{\rm if}\ p\in\mathbb{Z},\\[6pt]
2Y_{p+r},&{\rm if}\ p\in\frac{1}{2}+\mathbb{Z},\end{array}\right.
\ \,[G_r,M_p]=\left\{\begin{array}{cc}
\frac{p}{2}M_{p+r},&{\rm if}\ p\in\mathbb{Z},\\[6pt]
2M_{p+r},&{\rm if}\ p\in\frac{1}{2}+\mathbb{Z},\end{array}\right.\\
&&[Y_p,Y_q]=\left\{\begin{array}{ccc}
\frac{1}{2}(q-p)M_{q+p},&{\rm if}\ p,q\in\mathbb{Z},\\[6pt]
\frac{q}{2}M_{q+p},&{\rm if}\ p\in\mathbb{Z},q\in\frac{1}{2}+\mathbb{Z},\\[6pt]
2M_{q+p},&{\rm if}\ p,q\in\frac{1}{2}+\mathbb{Z}.\end{array}\right.
\end{eqnarray*}
It is easy to see that $\widetilde{\frak{tsns}}$ is $\mathbb{Z}_2$-graded with $\widetilde{\frak{tsns}}=\widetilde{\frak{tsns}}_{\bar{0}}\oplus\widetilde{\frak{tsns}}_{\bar{1}}$, where
\begin{eqnarray*}
\widetilde{\frak{tsns}}_{\bar{0}}\!\!\!&=&\!\!\!{\rm{Span}}_\mathbb{C}\{L_n,Y_n,M_n,\frak{c}\,|\,n\in \mathbb{Z}\},\\
\widetilde{\frak{tsns}}_{\bar{1}}\!\!\!&=&\!\!\!{\rm{Span}}_\mathbb{C}\{G_r,Y_r,M_r\,|\,r\in\frac{1}{2}+\mathbb{Z}\}.
\end{eqnarray*}
The Cartan subalgebra (exactly the maximal toral subalgebra) of $\widetilde{\frak{tsns}}$ is $\widetilde{\frak{h}}={\mathbb C}L_0\oplus{\mathbb C}M_0\oplus{\mathbb C}\frak{c}$. It should be noted that $\widetilde{\frak{tsns}}_{\bar{0}}$ is precisely the well-known twisted Schr\"{o}dinger-Virasoro Lie algebra and the subalgebra $\widetilde{\frak{ns}}$ spanned by
$\{L_n,G_r,\frak{c}\,|\,n\in{\mathbb{Z}},
\,r\in\frac{1}{2}+\mathbb{Z}\}$ is the $N\!=\!1$ Neveu-Schwarz algebra. Denote $\widetilde{\frak{I}}=\frak{I}\oplus\frak{c}$ with $\frak{I}={\rm{Span}}_\mathbb{C}\{Y_p,M_p\,|\,p\in\frac{1}{2}\mathbb{Z}\}$. It is easy to see that both $\widetilde{\frak{I}}$ and $\frak{I}$ are ideals of $\widetilde{\frak{tsns}}$.

An $\widetilde{\frak{h}}$-diagonalizable module over $\widetilde{\frak{tsns}}$ is usually called {\it weight module}. If all weight spaces of a weight $\widetilde{\frak{tsns}}$-moudule are finite-dimensional, the module is called a {\it Harish-Chandra module}. We introduce the following notations:
\begin{eqnarray*}
\widetilde{\frak{tsns}}_{+}\!\!\!&=&\!\!\!
{\rm{Span}}_\mathbb{C}\{L_n,Y_p,M_p,G_r\,|\,n\in\mathbb{Z}^+,\,p\in\frac{1}{2}\mathbb{Z},\,r\in\frac{1}{2}+\mathbb{Z},\,p,r>0\},\\
\widetilde{\frak{tsns}}_{-}\!\!\!&=&\!\!\!
{\rm{Span}}_\mathbb{C}\{L_n,Y_p,M_p,G_r\,|\,n\in\mathbb{Z}^-,\,p\in\frac{1}{2}\mathbb{Z},r\in\frac{1}{2}+\mathbb{Z},\,p,r<0\},\\
\widetilde{\frak{tsns}}_{0}\!\!\!&=&\!\!\!{\rm{Span}}_\mathbb{C}\{L_0,Y_0,M_0,\,\frak{c}\}.
\end{eqnarray*}
Then $\widetilde{\frak{tsns}}$ admits the following triangular decomposition:
\begin{eqnarray*}
\widetilde{\frak{tsns}}\!\!\!&=&\!\!\!\widetilde{\frak{tsns}}_{+}\oplus\widetilde{\frak{tsns}}_{0}\oplus\widetilde{\frak{tsns}}_{-}.
\end{eqnarray*}

For any irreducible weight $\widetilde{\frak{tsns}}$-module $V$, $L_0$, $M_0$ and $\frak{c}$ must act as some complex numbers on it. Furthermore, $V$ has the weight space decomposition
$V=\oplus_{\lambda\in\mathbb{C}}V_\lambda$, where $V_\lambda=\{v\in V\,|\,L_0v=\lambda v\}$ is called {\it a weight space} with weight $\lambda$.
Denote the set of {\it weights $\lambda$} of $V$ by
${\rm supp}(V):=\{\lambda\in\mathbb{C}\,|\,V_\lambda\neq 0\}$,
which is called the {\it support} of $V$. If $V$ is an
irreducible weight $\widetilde{\frak{tsns}}$-module, then there exists some $\lambda\in\mathbb{C}$
such that ${\rm supp}(V)\subseteq\lambda+\frac{1}{2}\mathbb{Z}$.

An irreducible weight module $V$ is called {\it a pointed module}
if there exists a weight $\lambda\in\mathbb{C}$ such that ${\rm dim\,}V_\lambda=1$.
The following natural problem was firstly referred in \cite{X--AMSNS1997}:
\begin{prob}\label{prob}
Is any irreducible pointed module over the Virasoro
algebra a Harish-chandra module?
\end{prob}

An irreducible weight module $V$ is called {\it a mixed module} if
there exist $\lambda\in\mathbb{C}$ and $k\in\mathbb{Z}^*$ such that ${\rm
dim\,}V_\lambda=\infty$ and ${\rm dim\,}V_{\lambda+k}<\infty$. The following
conjecture was given in \cite{M--MS2004}:
\begin{conj}\label{conj}
There are no irreducible mixed modules over the
Virasoro algebra.
\end{conj}

The positive answers to the above question and conjecture were given in \cite{MZ--LA2007}. Such question and conjecture have already been solved on the truncated Virasoro algebras, the $W$-algebra $W(2,2)$, the twisted Schr\"{o}dinger-Virasoro algebra, the twisted Heisenberg-Virasoro algebra and the Neveu-Schwarz algebra in \cite{GL--JMP2010,LGZ--JMP2008,LS--AMSE2009,SS--AMSE2007,ZX--CA2012}.
In this note, we also give the positive
answers to the above question and conjecture for the twisted $N\!=\!1$ Schr\"{o}dinger-Neveu-Schwarz algebra.
Our main result is the following:
\begin{theo}\label{theo1304120956}
For any irreducible weight $\widetilde{\frak{tsns}}$-module $V$,
if there exists some $\lambda\in\mathbb{C}$ such that ${\rm dim\,}V_\lambda=\infty$, then ${\rm supp}(V)=\lambda+\frac{1}{2}\mathbb{Z}$ and ${\rm dim\,}V_{\lambda+k}=\infty$ for every $k\in\frac{1}{2}\mathbb{Z}$. In other words, for any irreducible weight $\widetilde{\frak{tsns}}$-module $V$, the condition that there exists some $\lambda\in\mathbb{C}$ such that $0<{\rm
dim\,}V_\lambda<\infty$ implies that $V$ is a Harish-Chandra module.
Then there are no irreducible mixed $\widetilde{\frak{tsns}}$-modules.
\end{theo}

Throughout this paper, we respectively denote $\mathbb{Z}^*$, $\mathbb{Z}^+$, $\mathbb{Z}^-$, $\mathbb{Z}_+$ and $\mathbb{Z}_-$ the sets of the nonzero, positive, negative, nonnegative and non-positive integers.
\vspace*{18pt}

\centerline{\bf2\ \ Proof of Theorem \ref{theo1304120956}}
\setcounter{section}{2}\setcounter{theo}{0} \setcounter{equation}{0}

\vspace*{8pt}

We first recall a main result about the irreducible weight Neveu-Schwarz-modules given in \cite{ZX--CA2012}, which we cited here as the following lemma:
\begin{lemm}\label{LemmNSM}
Let $V$ be an irreducible weight $\widetilde{\frak{ns}}$-module.
Assume that there exists $\lambda\in\mathbb{C}$ such that ${\rm dim\,}V_\lambda=\infty$. Then ${\rm supp}(V)=\lambda+\frac{1}{2}\mathbb{Z}$ and for every $k\in\frac{1}{2}\mathbb{Z}$, we have ${\rm dim\,}V_{\lambda+k}=\infty$.
\end{lemm}

Since $\frak{I}={\rm{Span}}_\mathbb{C}\{Y_p,M_p\,|\,p\in\frac{1}{2}\mathbb{Z}\}$ is a nontrivial ideal of $\widetilde{\frak{tsns}}$, $\frak{I} V$ is a submodule of $V$ for any $\widetilde{\frak{tsns}}$-module $V$, which gives the following lemma (such lemma was not used in \cite{ZX--CA2012}).

\begin{lemm}\label{IM=0M}
For any irreducible weight $\widetilde{\frak{tsns}}$-module $V$, we have $\frak{I} V=0$ or $\frak{I} V=V$.
\end{lemm}

For any irreducible weight $\widetilde{\frak{tsns}}$-module $V$, if we can prove $\frak{I} V=0$, then $V$ will degenerate to be an irreducible weight $\widetilde{\frak{ns}}$-module. Then Theorem \ref{theo1304120956} follows from Lemma \ref{LemmNSM} in this case. Denote the universal enveloping algebra of $\widetilde{\frak{tsns}}$ by $U(\widetilde{\frak{tsns}})$. For any irreducible weight $\widetilde{\frak{tsns}}$-module $V$, if there exists some $0\neq v\in V$ such that $\frak{I} v=0$, then $U(\widetilde{\frak{tsns}})v=V$. And the following lemma follows:
\begin{lemm}\label{1304120930}
Assume that $V$ is an irreducible weight module over $\widetilde{\frak{tsns}}$. If there exists some $0\neq v\in V$ such that $\frak{I} v=0$. Then $\frak{I}$ acts trivially on $V$ and $V$ is simply an irreducible weight module over $\widetilde{\frak{ns}}$.
\end{lemm}

It is easily to see that $\{G_{\frac{3}{2}},\,G_{\frac{1}{2}},\,Y_{\frac{1}{2}},\,M_{\frac{1}{2}}\}$ and $\{G_{-\frac{3}{2}},\,G_{-\frac{1}{2}},\,Y_{-\frac{1}{2}},\,M_{-\frac{1}{2}}\}$ respectively generate $\widetilde{\frak{tsns}}_{+}$ and $\widetilde{\frak{tsns}}_{-}$. Then the following lemma follows from the fact that both
highest and lowest weight $\widetilde{\frak{tsns}}$-modules are Harish-Chandra modules (similar to the Virasoro algebra case investigated in \cite{MZ--LA2007}).
\begin{lemm}\label{wtsnsHLHCHM}
Assume that there exists $\mu\in\mathbb{C}$ and a non-zero
element $v\in V_\mu$, such that
\begin{eqnarray*}
G_{\frac{3}{2}}v=G_{\frac{1}{2}}v=Y_{\frac{1}{2}}v=M_{\frac{1}{2}}v=0\ \ \ \mbox{or}\ \ \ G_{-\frac{3}{2}}v=G_{-\frac{1}{2}}v=Y_{-\frac{1}{2}}v=M_{-\frac{1}{2}}v=0.
\end{eqnarray*}
Then $V$ is a Harish-Chandra module.
\end{lemm}

\noindent{\it Proof of Theorem \ref{theo1304120956}}\ \ \,We shall prove this theorem step by step by several lemmas.

Assume now that $V$ is an irreducible weight $\widetilde{\frak{tsns}}$-module such that there exists $\lambda\in\mathbb{C}$ satisfying ${\rm dim\,}V_\lambda=\infty$. Denote the set $S_\lambda\!=\!\{p\,|\,{\rm dim\,}V_{\lambda+p}<\infty,\,p\in\frac{1}{2}\mathbb{Z}^*\}$.
\begin{lemm}\label{1304052141}
For any $\lambda\in\mathbb{C}$ satisfying ${\rm dim\,}V_\lambda=\infty$, there are at most two adjacent elements in $S_\lambda$. For convenience, we can suppose $S_\lambda\subseteq\{\frac12,\,1\}$.
\end{lemm}
{\it Proof}\ \ \ Suppose there are two different elements $p$ and $q$ in $S_\lambda$. Without loss of generality, we can assume $p=\frac12$, $q>\frac12$.

\noindent{\bf Case 1\ \ }$p=\frac12$, $q>\frac12$ and $q\in\frac{1}{2}+\mathbb{Z}$.

Let $W$ be the intersection of the kernels of the linear maps $X_{\frac{1}{2}}:\,V_\lambda\rightarrow V_{\lambda+\frac{1}{2}}$ and $X_{q}:\,V_\lambda\rightarrow V_{\lambda+q}$ for $X=G,\,Y,\,M$. Then $X_{\frac{1}{2}}W=X_{q}W=0$ for all $X\in\{G,\,Y,\,M\}$. The assumptions ${\rm dim\,}V_\lambda=\infty$, ${\rm dim\,}V_{\lambda+\frac{1}{2}}<\infty$ and ${\rm dim\,}V_{\lambda+q}<\infty$ force ${\rm dim\,}W=\infty$. Using $G_{\frac{1}{2}}W=Y_{\frac{1}{2}}W=M_{\frac{1}{2}}W=0$, $G_{q}W=Y_{q}W=M_{q}W=0$ and the given brackets, we obtain $L_nW=G_rW=Y_pW=M_pW=0$ for all $n=1,2q,2q+1,2q+2,\cdots$, $r=\frac12,q,q+1,q+2,\cdots$, $p\in\frac12\mathbb{Z}_{>0}$. We claim that $q>\frac32$ and $G_\frac32w\neq0$ for all $0\neq w\in W$. Otherwise, $W$ would be a Harish-Chandra module by Lemma \ref{wtsnsHLHCHM}. Thus ${\rm dim\,}G_\frac32W=\infty$. Since ${\rm dim\,}V_{\lambda+\frac{1}{2}}<\infty$, there exists some $w\in W$ such that $v=G_\frac32w$ and $L_{-1}v=Y_{-1}v=M_{-1}v=0$. It is easily to verify that $\widetilde{\frak{tsns}}_+v=0$, which forces $V$ to be a Harish-Chandra module according to Lemma \ref{wtsnsHLHCHM}. This contracts with our assumptions. Hence this case is impossible.

\noindent{\bf Case 2\ \ }$p=\frac12$, $q>1$ and $q\in\mathbb{Z}$.

Let $W$ be the intersection of the kernels of the linear maps $X_{\frac{1}{2}}:\,V_\lambda\rightarrow V_{\lambda+\frac{1}{2}}$ and $Z_{q}:\,V_\lambda\rightarrow V_{\lambda+q}$ for $Z=L,\,Y,\,M$. Then $X_{\frac{1}{2}}W=Z_{q}W=0$ for all $X\in\{G,\,Y,\,M\}$ and $Z\in\{L,\,Y,\,M\}$. The assumptions ${\rm dim\,}V_\lambda=\infty$, ${\rm dim\,}V_{\lambda+\frac{1}{2}}<\infty$ and ${\rm dim\,}V_{\lambda+q}<\infty$ force ${\rm dim\,}W=\infty$. Using $G_{\frac{1}{2}}W=Y_{\frac{1}{2}}W=M_{\frac{1}{2}}W=0$, $L_{q}W=Y_{q}W=M_{q}W=0$ and the given brackets, we obtain $L_nW=G_rW=Y_pW=M_pW=0$ for all $n=1,q,q+1,q+2,\cdots$, $r=\frac12,q+\frac{1}{2},q+\frac{3}{2},\cdots$, $p\in\frac12\mathbb{Z}^+$. We claim that $G_\frac32w\neq0$ for all $0\neq w\in W$. Otherwise, $W$ would be a Harish-Chandra module by Lemma \ref{wtsnsHLHCHM}. Thus ${\rm dim\,}G_\frac32W=\infty$. Similar to the proof of the first case, one can prove that this case is also impossible. Then this lemma follows.\hfill$\Box$\par

The following lemma can be easily verified, partially given in Lemma 2 of \cite{ZX--CA2012}.
\begin{lemm}\label{vMGeq0}
(i)\ \ Let $0\neq v\in V$ be such that $G_{\frac{1}{2}}v=Y_{\frac{1}{2}}v=M_{\frac{1}{2}}v=0$. Then
\begin{eqnarray*}
&&\big(\frac{1}{2}L_1G_{\frac{1}{2}}-G_{\frac{3}{2}}\big)G_{\frac{3}{2}}v
=0=Y_pv=M_pv\ \ {\rm for\ all}\ p\in\frac12\mathbb{Z}^{+}.
\end{eqnarray*}
(ii)\ \ Let $0\neq v\in V$ be such that $G_{-\frac{1}{2}}v=Y_{-\frac{1}{2}}v=M_{-\frac{1}{2}}v=0$. Then
\begin{eqnarray*}
&&\big(\frac{1}{2}L_{-1}G_{-\frac{1}{2}}+G_{-\frac{3}{2}}\big)G_{-\frac{3}{2}}v
=0=Y_{-p}v=M_{-p}v\ \ {\rm for\ all}\ p\in\frac12\mathbb{Z}^{-}.
\end{eqnarray*}
\end{lemm}

According to Lemma \ref{1304052141}, we can assume $V$ is an irreducible weight $\widetilde{\frak{tsns}}$-module such that there exists some $\mu\in\mathbb{C}$ satisfying ${\rm dim\,}V_\mu<\infty$, ${\rm dim\,}V_{\mu+\frac12}<\infty$  and ${\rm dim\,}V_{\mu+p}=\infty$ for all $p\in\mathbb{Z}^*\cup\{\frac12+\mathbb{Z}^*\}$.

\begin{lemm}\label{mu-1frac12}
$\mu\in\{-1,\,\frac12\}$.
\end{lemm}
{\it Proof}\ \ \
Let $U$ be the intersection of the kernels of the linear maps $X_{\frac{1}{2}}:\,V_{\mu-\frac{1}{2}}\rightarrow V_{\mu}$ for $X=G,\,Y,\,M$. Then $G_{\frac{1}{2}}U=Y_{\frac{1}{2}}U=M_{\frac{1}{2}}U=0$. The assumptions ${\rm dim\,}V_\mu<\infty$ and ${\rm dim\,}V_{\mu-{\frac{1}{2}}}=\infty$ force ${\rm dim\,}U=\infty$. We claim that $G_{\frac{3}{2}}u\neq0$ for all $0\neq u\in U$. Otherwise, $M$ would be a Harish-Chandra module by Lemma \ref{wtsnsHLHCHM}. Thus ${\rm dim\,}G_{\frac{3}{2}}U=\infty$. Since ${\rm dim\,}M_{\mu}<\infty$, there exists some $u\in U$ such that $v=G_{\frac{3}{2}}u$ and $G_{-\frac{1}{2}}v=Y_{-\frac{1}{2}}v=M_{-\frac{1}{2}}v=0$. Using Lemma \ref{vMGeq0}, we obtain $(\frac{1}{2}L_1G_{\frac{1}{2}}-G_{\frac{3}{2}})v=0$ and $Y_pu=M_pu=0$ for all $p\in\frac12\mathbb{Z}^+$. Then $U$ is simply a $\widetilde{\frak{ns}}$-module. Then we can use the same discussions as those given in Lemma 3 of \cite{ZX--CA2012}, and have the following observations:
\begin{eqnarray*}
0=L_{-1}G_{-\frac{1}{2}}(\frac{1}{2}L_1G_{\frac{1}{2}}-G_{\frac{3}{2}})w=(2L^2_0-3L_0)w
=\big(2(\mu+1)^2-3(\mu+1)\big).
\end{eqnarray*}
Then the lemma follows.\hfill$\Box$\par

According to Lemma \ref{mu-1frac12}, any irreducible weight $\widetilde{\frak{tsns}}$-module $V$ possessing only two different finite-dimensional weight spaces can be assumed to be of the following type: ${\rm dim\,}V_{\frac12}<\infty$, ${\rm dim\,}V_{1}<\infty$ and ${\rm dim\,}V_{p+\frac12}=\infty$ for all $p\in\mathbb{Z}^*\cup\{\frac12+\mathbb{Z}^*\}$. Without loss of generality, we will just prove that this case can not happen, which can be formulated as the following lemma:
\begin{lemm}\label{frac121nex}
The irreducible weight $\widetilde{\frak{tsns}}$-module $V$ satisfying ${\rm dim\,}V_\frac12<\infty$, ${\rm dim\,}V_{1}<\infty$  and ${\rm dim\,}V_{\frac12+p}=\infty$ for all $p\in\mathbb{Z}^*\cup\{\frac12+\mathbb{Z}^*\}$ does not exist.
\end{lemm}
{\it Proof}\ \ \ Let $W$ be the intersection of the kernels of the linear maps $X_{\frac{1}{2}}:\,V_{0}\rightarrow V_{\frac{1}{2}}$ for all $X\in\{G,\,Y,\,M\}$. Then $G_{\frac{1}{2}}W=Y_{\frac{1}{2}}W=M_{\frac{1}{2}}W=L_{0}W=0$ and ${\rm dim\,}W=\infty$.
According to the following facts:
\begin{eqnarray*}
&&[G_{\frac12},G_{\frac12}]=2L_1,\ \ \,[L_1,Y_r]=rY_{1+r},
\ \ \,[L_1,M_r]=(r+\frac12)M_{r+1},\ \ \,\forall\,\,r\in\frac12+\mathbb{Z}_+,
\end{eqnarray*}
we can prove the following identities:
\begin{eqnarray}\label{1304120750}
&&L_1W=Y_rW=M_rW=0,\ \ \,\forall\,\,r\in\frac12+\mathbb{Z}_+.
\end{eqnarray}
Recalling the following identities:
\begin{eqnarray*}
&&[G_{\frac12},Y_p]=2Y_{p+\frac12},
\ \ \,[Y_p,Y_q]=2M_{p+q},\ \ \,\forall\,\,p,\,q\in\frac12+\mathbb{Z},
\end{eqnarray*}
we can obtain the following identities:
\begin{eqnarray}\label{1304120751}
&&Y_nW=M_nW=0,\ \ \,\forall\,\,n\in\mathbb{Z}_+.
\end{eqnarray}
Then combining \eqref{1304120750} and \eqref{1304120751}, we arrive at the following identities:
\begin{eqnarray}\label{1304120753}
&&Y_pW=M_pW=0,\ \ \,\forall\,\,p\in(\frac12\mathbb{Z}_+)^*.
\end{eqnarray}
Using the above obtained results, we can prove the following identities:
\begin{eqnarray}\label{1304120757}
&&Y_pG_{\frac32}W=M_pG_{\frac32}W=0,\ \ \,\forall\,\,p\in(\frac12\mathbb{Z}_+)^*.
\end{eqnarray}
We claim that $G_{\frac{3}{2}}w\neq0$ for all $0\neq w\in W$. Otherwise, $V$ would become a Harish-Chandra module by Lemma \ref{wtsnsHLHCHM}, which forces
${\rm dim\,}G_{\frac{3}{2}}W=\infty$. Recalling our suppose ${\rm dim\,}V_{1}<\infty$, we can deduce that there exists $0\neq v\in G_{\frac{3}{2}}W$ such that $G_{-\frac12}v=Y_{-\frac12}v=M_{-\frac12}v=0$. Then $L_{-1}v=G_{-\frac12}G_{-\frac12}v=0$. According to the following identities:
\begin{eqnarray*}
&[L_{-1},Y_p]=pY_{p-1},&[L_{-1},M_p]=(p-\frac12)M_{p-1},\\
&[G_{-\frac12},Y_p]=2Y_{p-\frac12},
&[Y_p,Y_q]=2M_{p+q},\ \ \,\forall\,\,p,\,q\in\frac12+\mathbb{Z}^-,
\end{eqnarray*}
we can deduce the following results:
\begin{eqnarray}\label{1304120921}
&&Y_pv=M_pv=0,\ \ \,\forall\,\,p,\,q\in(\frac12\mathbb{Z}_-)^*.
\end{eqnarray}
The identities given in \eqref{1304120757}, imply
\begin{eqnarray}\label{1304120922}
&&Y_pv=M_pv=0,\ \ \,\forall\,\,p,\,q\in(\frac12\mathbb{Z}_+)^*,
\end{eqnarray}
Combining $Y_0=\frac12[G_{-\frac12},Y_{\frac12}]$, $M_0=\frac12[G_{-\frac12},M_{\frac12}]$, \eqref{1304120921} and \eqref{1304120922},
we can deduce
\begin{eqnarray*}
&&Y_0v=M_0v=0,
\end{eqnarray*}
which together with \eqref{1304120921} and \eqref{1304120922}, gives
\begin{eqnarray*}
&&Y_pv=M_pv=0,\ \ \,\forall\,\,p\in\frac12\mathbb{Z}.
\end{eqnarray*}
Then recalling Lemma \ref{1304120930}, we can finally deduce that $Y_p, M_p$ act trivially on the whole $V$ for all $p\in\frac12\mathbb{Z}$ and $V$ is simply an irreducible weight module over $\widetilde{\frak{ns}}$. Then this lemma follows from Lemma \ref{LemmNSM}.\hfill$\Box$\par

By Lemmas \ref{1304052141}, \ref{mu-1frac12} and \ref{frac121nex}, we get the following lemma immediately.
\begin{lemm}\label{frac121nex1}
There is at most one element in $S_{\lambda}$.
\end{lemm}
According to Lemma \ref{frac121nex1}, we can suppose there exists some $p\in\frac12\mathbb{Z}^*$, such that ${\rm dim\,}V_{\lambda+p}<\infty$. For convenience, we denote $\lambda+p$ by $\mu$ in the following lemma.
\begin{lemm}\label{frac121nex2}
$S_{\lambda}=\emptyset$.
\end{lemm}
{\it Proof}\ \ \,According to our suppose, we know ${\rm dim\,}V_{\mu}<\infty$ with $\mu\in S_{\lambda}$. Let $W$ be the intersection of the kernels of the linear maps $X_{\frac{1}{2}}:\,V_{\mu-\frac12}\rightarrow V_{\mu}$ for $X\in\{G,\,Y,\,M\}$. Then $G_{\frac{1}{2}}W=Y_{\frac{1}{2}}W=M_{\frac{1}{2}}W=L_{0}W=0$ and ${\rm dim\,}W=\infty$.
Similar to the corresponding proof of Lemma \ref{frac121nex}, we can prove the following identities:
\begin{eqnarray}
&&{\rm dim\,}G_{\frac{3}{2}}W=\infty,\label{1304121000}\\
&&Y_pG_{\frac32}W=M_pG_{\frac32}W=0,
\ \ \,\forall\,\,p\in(\frac12\mathbb{Z}_+)^*.\label{1304121001}
\end{eqnarray}
Recalling our suppose that ${\rm dim\,}V_{\mu}<\infty$, we claim that there exists some nonzero element $v\in G_{\frac{3}{2}}W$ such that
\begin{eqnarray}\label{1304121002}
&&L_{-1}v=Y_{-1}v=M_{-1}v=0.
\end{eqnarray}
According to \eqref{1304121001}, we also know that
\begin{eqnarray}\label{1304121003}
&&Y_pv=M_pv=0,\ \ \,\forall\,\,p,\,q\in(\frac12\mathbb{Z}_+)^*.
\end{eqnarray}
The following identities hold:
\begin{eqnarray*}
[L_{-1},Y_{1-p}]\!\!\!&=&\!\!\!
\left\{\begin{array}{lll}
(\frac{3}{2}-p)Y_{-p},&{\rm if}\ \,p\in\mathbb{Z},\\[8pt]
(1-p)Y_{-p},&{\rm if}\ \,p\in\frac{1}{2}+\mathbb{Z},
\end{array}\right.
\end{eqnarray*}
from which we can deduce
\begin{eqnarray}\label{1304121005}
&&Y_{p}v=0,\ \ \,\forall\,\,p\in\frac12\mathbb{Z}_-.
\end{eqnarray}
Combining \eqref{1304121003} and \eqref{1304121005}, we arrive at the following result:
\begin{eqnarray}\label{1304121005}
&&Y_{p}v=0,\ \ \,\forall\,\,p\in\frac12\mathbb{Z}.
\end{eqnarray}
The identity $M_pv=0$ for all $p\in\frac12\mathbb{Z}$ follows from \ref{1304121005} and the follows identities:
\begin{eqnarray*}
[Y_p,Y_q]\!\!\!&=&\!\!\!
\left\{\begin{array}{lll}
\frac{q}{2}M_{q+p},&{\rm if}\ p\in\mathbb{Z},q\in\frac{1}{2}+\mathbb{Z},\\[8pt]
2M_{q+p},&{\rm if}\ p,q\in\frac{1}{2}+\mathbb{Z}.
\end{array}\right.
\end{eqnarray*}
Then again recalling Lemma \ref{1304120930}, we can finally deduce that $Y_p, M_p$ act trivially on the whole $V$ for all $p\in\frac12\mathbb{Z}$ and $V$ is simply an irreducible weight module over $\widetilde{\frak{ns}}$. Then this lemma follows from Lemma \ref{LemmNSM}.\hfill$\Box$\par
Then Theorem \ref{theo1304120956} follows immediately from Lemma \ref{frac121nex2}.
\hfill$\blacksquare$\par

\vskip8pt

\noindent{\bf Acknowledgements}\ \,Supported by a NSF grant BK20160403 of Jiangsu Province and NSF grants 11501417, 11671056, 11271056, 11101056 of China, Innovation Program of Shanghai Municipal Education Commission,  Program for Young Excellent Talents in Tongji University and the Fundamental Research Funds for the Central Universities.


\vskip8pt

\begin{thebibliography}{9999}\parskip0pt\lineskip4pt\small

\bibitem{GL--JMP2010} X. Guo, X. Liu, Weight modules with a finite-dimensional weight space over the truncated Virasoro algebras, {\it J. Math. Phys}, {\bf 51} (2010), 123522.

\bibitem{H--JSM1994} M. Henkel, Schr\"{o}dinger invariance and strongly anisotropic critical systems, {\it J. Stat. Phys.}, {\bf 75}(1994), 1023--1029.

\bibitem{HU--NPB2006} M. Henkel, J. Unterberger, Supersymmetric extensions of Schr\"{o}dinger-invariance, {\it Nucl. Phys. B}, {\bf 746} (2006), 155--201.

\bibitem{LGZ--JMP2008} D. Liu, S. Gao, L. Zhu, Classification of irreducible weight modules over $W$-algebra $W(2,2)$, {\it J. Math. Phys}, {\bf 49} (2008), 113503.

\bibitem{LS--AMSE2009} J. Li, Y. Su, Irreducible weight modules over the twisted Schr\"{o}dinger-Virasoro algebra, {\it Acta Mathematica Sinica, English Series}, {\bf 25}(4) (2009), 531--536.

\bibitem{M--MS2004} V. Mazorchuk, On simple mixed modules over the Virasoro algebra, {\it Mat. Stud.}, {\bf 22} (2004), 121--128.

\bibitem{MZ--LA2007} V. Mazorchuk, K. Zhao, Classification of simple weight Virasoro modules with a finite-dimensional weight space, {\it J. Algebra}, {\bf307} (2007), 209--214.

\bibitem{RU--AHP2006} J. Unterberger, C. Roger, The Schr\"{o}dinger-Virasoro Lie algebra, {\it Springer Texts and Monographs in Physics}, Springer (Heidelberg 2012).

\bibitem{SS--AMSE2007} R. Shen, Y. Su, Classification of irreducible weight modules with a finite-dimensional weight space over twisted Heisenberg-Virasoro algebra, {\it Acta Math. Sinica, English Series}, {\bf 23} (2007), 189--192.

\bibitem{X--AMSNS1997} X. Xu, Pointed representations of Virasoro algebra, A Chinese summary appears in {\it Acta Math. Sinica, Chinese Series}, {\bf 40}(3), 479; {\it Acta Math. Sinica, New Series}, {\bf 13} (1997), 161--168.

\bibitem{ZX--CA2012} X. Zhang, Z. Xia, Classification of simple weight modules for the Neveu-Schwarz algebra with a finite-dimensional weight space, {\it Comm. Alg.} {\bf 40} (2012), 2161--2170.


\end{thebibliography}
\end{document}